\documentclass{amsart}

\newtheorem{theorem}{Theorem}[section]
\newtheorem{lemma}[theorem]{Lemma}
\newtheorem{prop}[theorem]{Proposition}

\theoremstyle{definition}

\theoremstyle{remark}
\newtheorem{remark}[theorem]{Remark}

\numberwithin{equation}{section}

\usepackage{amscd}
\usepackage{xypic}  
\usepackage{amssymb}
\usepackage{amsthm}
\usepackage{epsf}
\usepackage[T1]{fontenc}
\def\i{^{-1}}

\def\ker{\operatorname{ker}}

\def\min{\operatorname{min}}
\def\log{\operatorname{log}}

\def\c1{\operatorname{c_1}}
\def\c2{\operatorname{c_2}}
\def\coker{\operatorname{coker}}

\def\Sym{\operatorname{Sym}}
\def\Hilb{\operatorname{Hilb}}

\def\CC{{\mathbb C}}

\def\PP{{\mathbb P}}
\def\L{{\mathcal L}}
\def\Q{{\mathcal Q}}
\def\A{{\mathcal A}}
\def\D{{\mathfrak D}}
\def\B{{\mathcal B}}

\def\N{{\mathcal N}}
\def\O{{\mathcal O}}
\def\I{{\mathcal J}}

\def\E{{\mathcal E}}
\def\C{{\mathcal C}}
\def\H{{\mathcal H}}
\def\F{{\mathcal F}}

\def\M{{\mathcal M}}

\def\T{{\mathcal T}}
\def\i{{\mathfrak I}}
\def\j{{\mathfrak J}}
\def\x{\times}                   
\def\iso{\simeq}
\def\eqv{\equiv}
\def\sub{\subseteq}

\def\+{\oplus}                   
\def\*{\otimes}                  
\def\hpil{\longrightarrow}       
\def\khpil{\rightarrow}

\def\mod{\operatorname{mod}}

\def\Ext{\operatorname{Ext}}

\def\Shom{\operatorname{ \mathcal{H}\mathit{om} }}
\def\Shext{\operatorname{ \mathcal{E}\mathit{xt} }}

\def\Pic{\operatorname{Pic}}

\def\Sing{\operatorname{Sing}}

\def\hs{\hspace{.05in}}

\begin{document}

\title[Isolated smooth curves in Calabi-Yau threefolds]{On isolated smooth curves of low genera in Calabi-Yau complete intersection threefolds} 

\author{Andreas Leopold Knutsen}
\address{Andreas Leopold Knutsen, Department of Mathematics, University of Bergen,
Johannes Brunsgate 12, 5008 Bergen, Norway.}
\email{andreas.knutsen@math.uib.no}




\begin{abstract}
Building on results of Clemens and Kley, we find 
criteria for a continuous family of curves in a nodal $K$-trivial threefold $Y_0$ to deform to a scheme of finitely many smooth isolated curves in a general deformation $Y_t$ of $Y_0$. As an application, 
we show the existence of smooth isolated curves of bounded genera and unbounded 
degrees in Calabi-Yau complete intersections threefolds. 
\end{abstract}

\maketitle

\section{Introduction} 
\label{sec:intro}

In this paper we study embeddings of complex projective curves
into (smooth) Calabi-Yau complete intersection ($CICY$) threefolds. Such embeddings,
and Calabi-Yau threefolds in general, have in the past decades been
objects of extended interest in both algebraic geometry and physics.
The goal of counting such curves (especially rational) has inspired
the development of quantum cohomology and led to the discovery of
surprising relations between algebraic geometry and the theory of
mirror symmetry.

In the paper \cite{kl1}, Kley developed a framework for showing existence of curves of certain genera and degrees in $CICY$ threefolds. The paper built on the original idea in the case of genus zero curves of Clemens \cite{cl1} (then used also in \cite{katz}, \cite{og}
and \cite{EJS}): one starts with a $K3$ complete intersection surface $X$ 
containing a smooth rational curve $C$, embeds the surface in a nodal $CICY$ of suitable intersection type $Y$ and proves that under a general deformation $Y_t$ of $Y_0=Y$, the rational curve deforms to an isolated curve in the deformation. In the higher genus case, the curve $C$ is replaced by a complete linear system $|\L|$ of curves on the surface of dimension equal to the genus, and the idea is to prove that only finitely many of these deform to the deformation $Y_t$ and possibly also that these are smooth and isolated. 
The main existence result in \cite{kl1} is \cite[Thm. 1]{kl1}, claiming that 
for any $d \geq 3$, the general $CICY$ threefold contains smooth, isolated elliptic curves of degree $d$, except for degree $3$ curves in the $CICY$ of type $(2,2,2,2)$.

A crucial point in this construction is to show that the curves on the $K3$ surface $X$ do not acquire any additional deformations when considered as curves in $Y$, precisely that 
\[
h^0(\N_{C/X}) = h^0(\N_{C/Y}) \; \mbox{for all} \; C \in |\L|. 
\]
Unfortunately, the proof of this  step, namely \cite[Thm. 3.5]{kl1}, contains a serious gap, which also influences the proof of its corollary \cite[Cor. 3.6]{kl1}, which in fact cannot hold (cf. Remark \ref{rem:hf} for a more detailed explanation). As a consequence, the proof of \cite[Thm. 1]{kl1} is incorrect.

This paper has two main purposes:
\begin{itemize}
\item We give criteria for a continuous family of curves on a regular surface in a nodal
threefold $Y$ with trivial canonical bundle to deform to a scheme of finitely many smooth isolated curves in a general deformation $Y_t$ of $Y_0=Y$, using results from \cite{ck} and ideas from the unpublished preprint \cite{kl2} of Kley, see Theorem \ref{holger}. 
\item We apply these results to prove existence of smooth, 
isolated curves of low genera in the various $CICY$ threefold types, see Theorem \ref{result} (of which \cite[Thm. 1]{kl1} is a special case). 
 \end{itemize}

The first main result is the following. It is an improvement under slightly stronger hypotheses of a result in the preprint 
\cite{kl2} of Kley, which has never been published, presumably because of the gap in \cite{kl1}.

We first state the assumptions.
\vspace{0.4cm}

{\bf Setting and assumptions.}
Let $P$ be a smooth projective variety of dimension $r \geq 4$ and
$\E$ a vector bundle of
rank $r-3$ on $P$ that splits as a direct sum of line bundles
\[ \E = \+_{i=1}^{r-3} \M_i. \]
Let
\[ 
s_0=s_{0,1} \+ \cdots \+ s_{0,r-3} \in H^0 (P, \E)= \+_{i=1}^{r-3} H^0(P,\M_i) 
 \]
be a regular section, where $s_{0,i} \in H^0(P,\M_i)$ for $i=1,\ldots,r-3$. 
Set
\[ Y = Z (s_0) \; \; \mbox{and}  \; \; Z=Z(s_{0,1} \+ \cdots \+ s_{0,r-4}) \]
(where $Z=P$ if $r=4$).

Let $X \subset Y$ be a smooth, regular surface (i.e. $H^1(X, \O_X)=0$) and 
$\L$ a line bundle on $X$.

We make the following additional assumptions:

\begin{itemize}
\item [(A1)] $Y$ has trivial canonical bundle;
\item [(A2)] $Z$ is smooth along $X$ and the only singularities of $Y$ which lie in $X$ are $\ell$
  nodes $\xi_1, \ldots, \xi_{\ell}$.  Furthermore
  \[ \ell \geq \dim |\L| +2; \]
\item [(A3)] $|\L| \neq \emptyset$ and the general element of $|\L|$ is a smooth, irreducible curve; 
\item [(A4)] for every $\xi_i \in S:=\{\xi_1, \ldots, \xi_{\ell}\}$, if 
$|\L \* \I_{\xi_i}| \neq \emptyset$, then its general member is nonsingular at 
$\xi_i$; 
\item [(A5)] $H^0(C, \N _{C/X}) \iso H^0(C, \N _{C/Y})$ for all
            $C \in |\L|$;
\item [(A6)] $H^1(C, \N _{C/P})=0$ for all $C \in |\L|$;
\item [(A7)] the image of the natural restriction map 
\[  
\xymatrix{
H^0(P, \M_{r-3}) \ar[r] & H^0(S, \M_{r-3} \* \O_S) \iso \CC^{\ell} 
}
\]
has codimension one.
\end{itemize}

Let $s \in H^0(P,\E)$ be a general section.
Then our result is the following:

\begin{theorem} \label{holger}
  Under the above setting and assumptions (A1)-(A7), the members of
  $|\L|$ deform to a length $\ell-2 \choose \dim |L|$ scheme of curves
  that are smooth and isolated in the general deformation $Y_t=Z(s_0+ts)$ of $Y_0=Y$. In particular, $Y_t$ contains a smooth, isolated curve that is a deformation of a curve in $|\L|$. 
\end{theorem}

This result improves \cite[Thm. 1.1]{kl2}, since the curves in that theorem are
not claimed to be smooth or isolated.

Our main application is Theorem \ref{result} right below, of which \cite[Thm. 1]{kl1} is the special case with $g=1$. Thus we give a correct proof of \cite[Thm. 1]{kl1} and, at the same time, we extend the result to genera $>1$.

\begin{theorem} \label{result}
 Let $d \geq 1$ and $g \geq 0$ be integers. Then in any of the
 following cases the general Calabi-Yau complete intersection threefold $Y$ of 
the given
 type contains an isolated, smooth curve of degree $d$ and genus $g$:
\begin{itemize}
\item[(a)] $Y=(5) \sub \PP^4$: $g=0$ and $d >0$; $g=1$ and $d \geq 3$;
$2 \leq g \leq 6$ and $d \geq g+3$; $7 \leq g \leq 9$ and $d \geq g+2$;
$g=10$ and $d \geq 11$; $11 \leq g \leq 22$ and $d \geq \frac{g+13}{2}$.
\item[(b)] $Y=(4,2) \sub \PP^5$: $g=0$ and $d >0$; $g=1$ and $d \geq 3$;
$g=2$ and $d \geq 5$; $3 \leq g \leq 8$ and $d \geq g+4$;
$9 \leq g \leq 11$ and $d \geq g+3$; 
$12 \leq g \leq 15$ and $d \geq \frac{g+16}{2}$. 
\item[(c)] $Y=(3,3) \sub \PP^5$: $g=0$ and $d >0$; $g=1$ and $d \geq 3$;
$g=2$ and $d \geq 5$; $3 \leq g \leq 7$ and $d \geq g+4$.
\item[(d)] $Y=(3,2,2) \sub \PP^6$: $g=0$ and $d >0$; $g=1$ and $d \geq 3$;
$g=2$ and $d \geq 5$; $g=3$ and $d \geq 7$; 
$4 \leq g \leq 10$ and $d \geq g+5$.
\item[(e)] $Y=(2,2,2,2) \sub \PP^7$: $g=0$ and $d >0$; $g=1$ and $d \geq 4$;
$g=2$ and $d \geq 6$; $g=3$ and $d \geq 7$.
\end{itemize}
\end{theorem}

We remark that the genus zero case of the theorem is already known by
\cite{katz,og,EJS}. In \cite{kl2} an existence result similar to Theorem \ref{result} was claimed, but only for {\it geometrically rigid, connected curves} (not necessarily smooth and isolated). But the proof of that result also relied on \cite[Thm. 3.5]{kl1}.

The paper is  divided into two parts in a natural way:
 
The first part, consisting of Sections \ref{sec:local}-\ref{sec:mov}, is devoted to the proof of Theorem \ref{holger}. In Sections \ref{sec:local}-\ref{sec:cy}
we 
study deformations of curves in a complete linear system $|\L|$ lying on a smooth surface $X$ embedded in a nodal threefold $Y$ with  emphasis on the cases of regular surfaces  in 
threefolds with trivial canonical bundle. Special attention is devoted to studying if the curves in $|\L|$ on $X$ 
acquire additional deformations when embedded in $Y$, that is, to studying condition (A5). The crucial result is Proposition \ref{lemma:r2}, which 
states that condition (A5) is equivalent to the condition
\begin{itemize}
\item [(A5)'] The set of nodes 
$S$ imposes independent conditions on $|\L|$, and the \linebreak natural 
map 
$\xymatrix{ \gamma_C:  \; H^0(C,\N_{X/Y}\* \O_C)  \ar[r] & H^1(C,\N_{C/X})}$
(cf. \eqref{eq:r14}) is an  isomorphism for all $C \in |\L|$. 
\end{itemize}
First of all, the conditions in (A5)' may be easier to check than condition (A5) directly. More importantly, however, the first of the two conditions in (A5)'
assures that the locus of curves in $|\L|$ passing through at least one node is a simple normal crossing (SNC) divisor (consisting of $\ell$ hyperplanes). This enables us to identify a certain sheaf $\Q$ of obstructions to deformation as the locally free sheaf of differentials with logarithmic poles along an SNC divisor, cf. \eqref{eq:excess}, \eqref{eq:snc} and \eqref{eq:kons2}. This is a crucial point to assure that  {\it smooth and isolated} curves survive in a general deformation $Y_t$ of $Y_0=Y$.

The proof of Theorem \ref{holger} is finished in Section \ref{sec:mov}, following the proof of \cite[Thm. 1.1]{kl2}.

In the second part, consisting of Sections \ref{sec:set} and \ref{sec:con}, we
apply Theorem \ref{holger} to the case of $K3$ surfaces in complete intersection Calabi-Yau threefolds to prove Theorem \ref{result}. For each of the complete intersection types in Theorem \ref{result}, there is a standard construction allowing to embed a $K3$ surface of one (or more) of the three complete intersection types $(4)$ in $\PP^3$, $(2,3)$ in $\PP^4$ and $(2,2,2)$ in $\PP^5$ into a {\it nodal}
$CICY$ threefold. This will be recalled in Section \ref{sec:set}. We are then in the setting of Theorem \ref{holger} with $X$ the $K3$ surface, $Y$ the $CICY$, $P$ a projective space and $\E$ the vector bundle corresponding to the complete intersection type of $Y$. 
All various complete intersection types and possible constructions are summarized in  Table \ref{tabella} in Section \ref{sec:set}.

By construction and Bertini's theorem, condition (A1) and the first part of condition (A2) are satisfied. In each of the cases in Theorem \ref{result}, we will then need to verify the remaining conditions (A2)-(A7), and this is the reason for the various numerical conditions on $d$ and $g$ in the theorem. To help the reader navigate through the proof, we now briefly explain how it works.

The existence of smooth curves of certain degrees and genera on the three types of complete intersection $K3$ surfaces is given by the existence results in 
\cite{mori} and \cite{kn1}, cf. Theorem \ref{thm:bncurves}. We set $\L$ to be the line bundle defined by the curves and an important point is that the existence results guarantee that $\L$ and $\O_X(1)$ are independent in $\Pic X$. 

Condition (A3) is  automatically satisfied, as well as condition (A4), by standard
arguments, cf. Lemma \ref{lemma:curvesonk3}. Condition (A7) is also 
satisfied by construction, cf. Lemma \ref{lemma:a7}. To check (A5) we will check the two conditions in (A5)'. We prove that the second one is 
satisfied in Proposition \ref{thm:gap}, and it is here that we need to use the fact that $\L$ and $\O_X(1)$ are independent in $\Pic X$.

Therefore, at the end, the conditions that will be responsible for the numerical constraints in Theorem \ref{result} are conditions (A2) and (A6), as well as the first condition in (A5)', namely that the set of nodes $S$ imposes independent conditions on the linear system $|\L|$. This is perhaps the most tricky condition to check, and we give
sufficient conditions for this to hold in 
Lemma \ref{lemma:tappato}. 

The numerical conditions we end up with are given in Proposition \ref{prop:checkcond} (in addition to the conditions in the existence result Theorem \ref{thm:bncurves}).
Finally,  a case-by-case study of these conditions will give Theorem \ref{result}.

\subsection*{Conventions and definitions}
The ground field is the field of complex numbers. We say a curve $C$
in a variety $V$ is 
{\it geometrically rigid} in $V$ if the space of
embedded deformations of $C$ in $V$ is zero-dimensional. If
furthermore this space is reduced, we say that $C$ is {\it isolated} or
{\it infinitesimally rigid} in $V$. From the infinitesimal
study of the Hilbert scheme of $V$, it follows that $C$ is
infinitesimally rigid if and only if $H^0(C, \N _{C/V})=0$.

A $K3$ surface is a smooth projective (reduced and irreducible) surface $X$ with trivial
canonical bundle and such that $H^1 (\O_X)=0$. In particular $h^2 (\O_X)= 1$
and $\chi (\O_X)= 2$.

A Calabi-Yau threefold $Y$ is a projective variety of dimension $3$ with
trivial canonical bundle and $h^1(\O_Y)=h^2(\O_Y)=0$. In this paper a
Calabi-Yau threefold will be at worst nodal.

\subsection*{Acknowledgements}
I thank T.~Johnsen, H.~Clemens, H.~P.~Kley, A.~F.~Lopez, S.~A.  Str{\o}mme, 
S.~Lekaus, R.~Ile and G.~Fl{\o}ystad for useful conversations on this subject.

I was made aware of the gap in \cite[Thm. 3.5]{kl1} (and a similar one in \cite[Example 4.3]{ck}, cf. Remark \ref{rem:hh}) by a referee of an earlier version of my preprint \cite{kn2} from 2001, where I applied the results of \cite{kl1,kl2} to Calabi-Yau threefolds that are complete intersections in certain homogeneous spaces. (As a consequence, \cite{kn2} has still not been published.) I must 
therefore thank this referee for having discovered the gap in \cite{kl1} and for a very carefully written referee report where he or she explained this in detail.

Finally, I must thank the referee of {\it this} paper for a very careful reading of the manuscript and many suggestions. In particular, he or she caught a mistake in the first version and suggested a rewrite 
resulting in the present Sections \ref{sec:local}-\ref{sec:cy}, 
making the paper, I believe, much clearer and easier to read.

\section{Curves through nodes on threefolds, local theory} 
\label{sec:local}

Let 
\[ Y = \{ (x,y,z,w) \in \CC^4 \; | \; xw-yz=0 \} \]
be (the analytic germ of) a nodal threefold singularity in affine $4$-space containing the plane
\[ X=\{(x,y,z,w) \in \CC^4 \; | \; x=y=0\}. \]

Let $\I_{X/Y} \subset \O_Y$ denote the ideal sheaf of $X$ in $Y$. A resolution of the $\O_X$-module $\I_{X/Y}$ is cyclic of the form

\begin{equation} \label{eq:cyclic}
\xymatrix{
& \cdots \hspace{-2cm} & \hspace{0,5cm} \O_Y^{\+ 2} \ar[r]^{ \left( \begin{array}{cc}
x  & z  \\
y & w \end{array} \right) } \hspace{0,5cm} 
 & \hspace{0,5cm} \O_Y^{\+ 2}  
\ar[r]^{ \left( \begin{array}{cc}
w  & -z  \\
-y & x \end{array} \right) } \hspace{0,5cm} &   
&  & 
\\
& &  \hspace{0,5cm}\O_Y^{\+ 2}  \ar[r]^{ \left( \begin{array}{cc}
x  & z  \\
y & w \end{array} \right) } \hspace{0,5cm} 
 & \hspace{0,5cm} \O_Y^{\+ 2}  
\ar[r]^{ \left( \begin{array}{cc}
w  & -z  \\
-y & x \end{array} \right) } \hspace{0,5cm} & \hspace{0,5cm} 
\O_Y^{\+ 2} \ar[r] ^{ \left( \begin{array}{c}
x  \\
y  \end{array} \right) } \hspace{0,5cm} 
& \hspace{0,1cm} \I_{X/Y} \ar[r] & 0
}
\end{equation}
so that we have 
\begin{equation}
  \label{eq:r1}
\Shext^i_{\O_Y}(\I_{X/Y}, \O_X) =  
\begin{cases} 
0, & \;  \mbox{if $i$ is even, } \\
\frac{\O_X}{(z,w)} \iso \CC, & \;  \mbox{if $i$ is odd.}
\end{cases} 
\end{equation}
Tensoring the sequence \eqref{eq:cyclic} with $\O_X$ is terminally exact and so we obtain the resolution
\begin{equation} \label{eq:res}
\xymatrix{
0 \ar[r] & \O_X \ar[r]^{ \left( \begin{array}{cc}
w  & -z  \end{array} \right) } \hspace{0,5cm} & \hspace{0,5cm} 
\O_X^{\+ 2} \ar[r] ^{ \left( \begin{array}{c}
x  \\
y  \end{array} \right) } \hspace{0,5cm} 
& \hspace{0,5cm} \frac{\I_{X/Y}}{\I_{X/Y}^2} \ar[r] & 0
}
\end{equation}
of $\frac{\I_{X/Y}}{\I_{X/Y}^2}$ as an $\O_X$-module. Thus
\[ \N_{X/Y} := \Shom_{\O_X} \Big(\frac{\I_{X/Y}}{\I_{X/Y}^2}, \O_X\Big) = \Big\{(a,b) \in \O_X^{\+2} \; | \; wa=zb\Big\} = \Big\{(cz,cw) \; | \; c \in \O_X \Big\} \]
is locally free. Moreover
\begin{equation}
  \label{eq:nodi1l}
 \Shext^1_{\O_X} \Big(\frac{\I_{X/Y}}{\I_{X/Y}^2}, \O_X\Big) = \frac{\O_X}{(z,w)} \iso \CC. 
\end{equation}

Let now
\[ C = \{(x,y,z,w) \in \CC^4 \; | \; x=y=f(z,w)=0\} \]
be a curve in $X$ passing through the node $(0,0,0,0)$ of $Y$. 

Tensoring \eqref{eq:res} by $\O_C$ we obtain the resolution
\begin{equation} \label{eq:resc}
\xymatrix{
0 \ar[r] & \O_C \ar[r]^{ \left( \begin{array}{cc}
w  & -z  \end{array} \right) } \hspace{0,5cm} & \hspace{0,5cm} 
\O_C^{\+ 2} \ar[r] ^{ \left( \begin{array}{c}
x  \\
y  \end{array} \right) } \hspace{0,5cm} 
& \hspace{0,5cm} \frac{\I_{X/Y}}{\I_{X/Y}^2} \* \O_C\ar[r] & 0
}
\end{equation}
of $\frac{\I_{X/Y}}{\I_{X/Y}^2} \* \O_C$ as an $\O_C$-module, as the map
$\left( \begin{array}{cc}
w  & -z  \end{array} \right)$ is injective.
Thus the map from \eqref{eq:res} to \eqref{eq:resc}
gives the commutative diagram
\begin{equation} \label{eq:r2}
\begin{footnotesize} \xymatrix{
& 0  \ar[d] & 0 \ar[d] & 0 \ar[d] & \\
0 \ar[r] & \O_X \ar[r]^{\cdot f} \ar[d] & \O_X \ar[r] \ar[d]& \O_C \ar[r] \ar[d]& 0 \\
0 \ar[r] & \O_X^{\+2} \ar[r]^{\cdot f} \ar[d] & \O_X^{\+2} \ar[r] \ar[d]& \O_C^{\+2} \ar[r] \ar[d] & 0  \\
0 \ar[r] & \frac{\I_{X/Y}}{\I_{X/Y}^2} \ar[r]^{\cdot f} \ar[d] & \frac{\I_{X/Y}}{\I_{X/Y}^2} \ar[r] \ar[d]& \frac{\I_{X/Y}}{\I_{X/Y}^2} \* \O_C \ar[r] \ar[d] & 0  \\
& 0   & 0  & 0  & 
}
\end{footnotesize}
\end{equation}
in which the first two columns are projective $\O_X$-resolutions of 
$\frac{\I_{X/Y}}{\I_{X/Y}^2}$ and the third is a projective $\O_C$-resolution of 
$\frac{\I_{X/Y}}{\I_{X/Y}^2} \* \O_C$. 

For all $i \geq 0$, we define
\begin{equation}
  \label{eq:defFiC}
  \F^i_C : = \Shext^i_{\O_C} \Big(\frac{\I_{X/Y}}{\I_{X/Y}^2} \* \O_C,\O_C\Big).
\end{equation}
Applying $\Shom_{\O_X}(-,\O_X)$ to the first two columns of \eqref{eq:r2}, we obtain
\[
\begin{footnotesize} \xymatrix{
0 \ar[r] & \Shom_{\O_X}(\O_X,\O_X) \ar[r]^{\cdot f}  & \Shom_{\O_X}(\O_X,\O_X) \ar[r]  & \Shom_{\O_C}(\O_C,\O_C) \ar[r]  & 0 \\
0 \ar[r] & \Shom_{\O_X}(\O_X^{\+2},\O_X) \ar[r]^{\cdot f} \ar[u] & \Shom_{\O_X}(\O_X^{\+2},\O_X) \ar[r] \ar[u] & \Shom_{\O_C}(\O_C^{\+2},\O_C) \ar[r] \ar[u] & 0  \\
0 \ar[r] & \Shom_{\O_X}(\frac{\I_{X/Y}}{\I_{X/Y}^2},\O_X) \ar[r]^{\cdot f} \ar[u] & 
\Shom_{\O_X}(\frac{\I_{X/Y}}{\I_{X/Y}^2},\O_X) \ar[r] \ar[u] & \F^0_C
\ar[u]   \\
& 0 \ar[u]   & 0 \ar[u]  & 0 \ar[u]   & 
}
\end{footnotesize}
\]
From the snake lemma  we obtain the  exact sequence
\[
\begin{footnotesize} \xymatrix{
0 \ar[r] & \Shom_{\O_X}(\frac{\I_{X/Y}}{\I_{X/Y}^2},\O_X) \ar[r]^{\cdot f}  & 
\Shom_{\O_X}(\frac{\I_{X/Y}}{\I_{X/Y}^2},\O_X) \ar[r]  & \F^0_C & \\
\ar[r] & \Shext^1_{\O_X}(\frac{\I_{X/Y}}{\I_{X/Y}^2},\O_X) \ar[r]^{\cdot f}  & 
\Shext^1_{\O_X}(\frac{\I_{X/Y}}{\I_{X/Y}^2},\O_X) \ar[r]  & \F^1_C  \ar[r] & 0.   
}
\end{footnotesize}
\]
Since $f$ is the local equation of $C$, we can, by tensoring with $\O_C$, deduce the 
short exact sequence
\begin{equation}
  \label{eq:r5}
\xymatrix{
0  \ar[r] & \N_{X/Y} \* \O_C \ar[r] & 
\F^0_C  \ar[r] & \Shext^1_{\O_X} \Big(\frac{\I_{X/Y}}{\I_{X/Y}^2} ,\O_X\Big) \*\O_C  \ar[r] & 0
}
\end{equation}
and the isomorphism
\begin{equation}
  \label{eq:r4}
\xymatrix{ \Shext^1_{\O_X} \Big(\frac{\I_{X/Y}}{\I_{X/Y}^2} ,\O_X\Big) \*\O_C 
\ar[r]^{\hspace{1.8cm}\iso} & \F^1_C. 
}
\end{equation}

\section{Curves through nodes on threefolds, global theory} 
\label{sec:global}

In this section $X$ will be a smooth projective surface, $Y$ a projective threefold and $P$ a smooth projective variety of dimension $\geq 4$ such that 
$X \subset Y \subset P$. We assume that the 
only singularities of $Y$ lying on $X$ are finitely many nodal singularities
and that the embedding $Y \subset P$ is regular (e.g., $Y$ is a complete intersection in $P$). We denote the set of nodes of $Y$ lying on $X$ by $S$. Note that we have
\begin{equation}
  \label{eq:nodi1}
 \Shext^1_{\O_X} \Big(\frac{\I_{X/Y}}{\I_{X/Y}^2}, \O_X\Big) \iso \O_S, 
\end{equation}
by \eqref{eq:nodi1l}.

 Let $\L$ be a line bundle on $X$ such that $|\L| \neq \emptyset$ and let $C \in |\L|$. We define the sheaves $\F^i_C$ as in \eqref{eq:defFiC}.

Applying $\Shom_{\O_C}(-,\O_C)$ to
\begin{equation}
  \label{eq:r6}
\xymatrix{
0  \ar[r] & \frac{\I_{X/Y}}{\I_{X/Y}\I_{C/Y}} \ar[r] & \frac{\I_{C/Y}}{\I_{C/Y}^2}
\ar[r] & \frac{\I_{C/X}}{\I_{C/X}^2} \ar[r] & 0
}
\end{equation}
we obtain the exact sequence
\begin{equation} \label{eq:diaF}
\xymatrix{
0 \ar[r] & \N_{C/X} \ar[r] & \N_{C/Y} \ar[r] & \F^0_C \ar[r] & 0
}
\end{equation}
as well as
\begin{equation} \label{eq:diaF2}
\Shext^1_{\O_C}\Big(\frac{\I_{C/Y}}{\I_{C/Y}^2},\O_C\Big) \iso 
\Shext^1_{\O_C}\Big(\frac{\I_{X/Y}}{\I_{X/Y}^2} \*\O_C,\O_C\Big) = \F^1_C. 
\end{equation}
Applying $\Shom_{\O_C}(-,\O_C)$ to
\begin{equation}
  \label{eq:r7}
\xymatrix{
0  \ar[r] & \frac{\I_{Y/P}}{\I_{X/P}\I_{C/P}} \ar[r] & \frac{\I_{X/P}}{\I_{X/P}\I_{C/P}} \ar[r] & \frac{\I_{X/Y}}{\I_{X/Y}\I_{C/Y}} \ar[r] & 0
}
\end{equation}
we obtain the exact sequence
\begin{equation} \label{eq:r9}
\xymatrix{
0 \ar[r] & \F^0_C \ar[r] & \N_{X/P} \* \O_C \ar[r] & \N_{Y/P} \*\O_C \ar[r] & 
\F^1_C \ar[r] & 0
}
\end{equation}
and applying  $\Shom_{\O_C}(-,\O_C)$ to
\begin{equation}
  \label{eq:r8}
\xymatrix{
0  \ar[r] & \frac{\I_{Y/P}}{\I_{Y/P}\I_{C/P}} \ar[r] & \frac{\I_{C/P}}{\I_{C/P}^2} 
\ar[r] & \frac{\I_{C/Y}}{\I_{C/Y}^2} \ar[r] & 0
}
\end{equation}
we obtain the exact sequence
\begin{equation} \label{eq:nod7}
\xymatrix{
0 \ar[r] & \N_{C/Y} \ar[r] & \N_{C/P} \ar[r] & \N_{Y/P} \* \O_C \ar[r] & 
\Shext^1_{\O_C} \Big(\frac{\I_{C/Y}}{\I_{C/Y}^2}, \O_C\Big) \ar[r] & 0,
}
\end{equation}
because $\Shext^1_{\O_C} \Big(\frac{\I_{C/P}}{\I_{C/P}^2}, \O_C\Big)=0$, as 
$C \subset P$ is a regular embedding.

Similarly, we also have the standard short exact sequence
\begin{equation} \label{eq:stand}
\xymatrix{
0 \ar[r] &  \N_{C/X} \ar[r] & \N_{C/P} \ar[r] & \N_{X/P} \* \O_C \ar[r] &  0.
} 
\end{equation}

Combining the sequences \eqref{eq:diaF}, \eqref{eq:r9}, \eqref{eq:nod7} and \eqref{eq:stand},
together with the isomorphisms  \eqref{eq:nodi1} and \eqref{eq:diaF2}, we obtain

\begin{equation} \label{eq:nod0r} 
\begin{footnotesize}\xymatrix{
         & 0   \ar[d]     & 0 \ar[d]  &  & &   \\
&         \N_{C/X} \ar[d] \ar@{=}[r] & \N_{C/X} \ar[d] & & &   \\
0 \ar[r] & \N_{C/Y} \ar[d] \ar[r] & \N_{C/P} \ar[d] \ar[r] 
&   \N_{Y/P} \* \O_C \ar@{=}[d]  \ar[r] &  \Shext^1_{\O_C} \Big(\frac{\I_{C/Y}}{\I_{C/Y}^2}, \O_C\Big)  \ar[r]  \ar[d]^{\cong} \ar[r] & 0  \\     
0 \ar[r] & \F^0_C  \ar[d] \ar[r] & \N_{X/P} \* \O_C  \ar[d] \ar[r]^{\Phi_C} 
&   \N_{Y/P} \* \O_C  \ar[r] &  \F^1_C  \iso \O_{S \cap C} \ar[r]  & 0  \\ 
         & 0       & 0   &  & &   
}
\end{footnotesize}
\end{equation}

\section{Curves on regular surfaces in $K$-trivial threefolds with nodes} 
\label{sec:cy}

In this section we keep the hypotheses and notation from the previous section and assume further that
\[ \omega_Y \iso \O_Y \; \; \mbox{and} \; \; h^1(\O_X)=0. \]
By adjunction
\begin{equation}
  \label{eq:r111}
\N_{X/Y} \iso \omega_X
\end{equation} 
except at the finite set $S$. But both sides are locally free so the isomorphism holds everywhere. Therefore, again by adjunction
\begin{equation}
  \label{eq:r110}
 \N_{X/Y} \* \O_C \iso \omega_C \* \N_{C/X}^{\vee}, 
\end{equation}
so that \eqref{eq:r5} reads
\begin{equation}
  \label{eq:r50}
\xymatrix{
0  \ar[r] & \omega_C \* \N_{C/X}^{\vee} \ar[r] & 
\F^0_C  \ar[r] & \O_{S \cap C}  \ar[r] & 0,
}
\end{equation}
where we have used \eqref{eq:nodi1}.

The main aim of this section is to find criteria for the curves $C \in |\L|$ not to obtain any additional deformations in $Y$, that is, for the equality $h^0(\N_{C/X})=h^0(\N_{C/Y})$ to hold for all $C \in |\L|$. 

We will need to consider, for $C \in |\L|$, the composition
\begin{equation}
  \label{eq:r14}
 \xymatrix{ \gamma_C: \;  H^0(C,\N_{X/Y}\* \O_C)
 \ar@{^{(}->}[r]^{\hspace{1.1cm}\iota_C}  &  H^0(C,\F^0_C)  
\ar[r] & H^1(C,\N_{C/X}),}  
\end{equation}
where the left hand inclusion $\iota_C$ arises from \eqref{eq:r5} (or, equivalently, \eqref{eq:r50})
 and the right hand map 
is the connecting homomorphism of \eqref{eq:diaF}. Note that
$H^0(C,\N_{X/Y}\* \O_C)$ and $H^1(C,\N_{C/X})$ are in fact dual by 
\eqref{eq:r110}. Therefore, $\gamma_C$ is an isomorphism if and only if it is of maximal rank.

\begin{lemma} \label{lemma:r1}
  Assume that $\omega_Y \iso \O_Y$ and $h^1(\O_X)=0$. 

The inclusion $\iota_C$ is an isomorphism for all $C \in |\L|$ if and only if 
$S$ imposes independent conditions on $|\L|$.
\end{lemma}

\begin{remark} \label{rem:snc}
  {\rm The condition that $S$ imposes independent conditions on $|\L|$ means, precisely, that if $\ell$ is the number of nodes, then $|\L \* \I_S| = \emptyset$ if $\ell > \dim |\L|$, and $\dim |\L \* \I_S| = \dim |\L|-\ell$ if $\ell \leq \dim |\L|$. This can also be rephrased as 
$\dim  |\L \* \I_{\xi_1} \* \I_{\xi_2} \* \cdots \* \I_{\xi_k}| = \dim |\L|-k$ for any subset $\{\xi_1, \ldots, \xi_k\}$ of $k$ nodes of $S$, whenever 
$k \leq \dim  |\L|$, or, even simpler, that at most $\dim |\L|$ of the nodes can lie on an element of $|\L|$.

Note that if $S$ imposes independent conditions on $|\L|$, then, in particular, the points in $S$ are different from the possible base points of $|\L|$, so that
the locus of curves in $|\L|$ passing through at least one node is an effective divisor
in $|\L|$ consisting of hyperplanes. (If $\dim |\L|=0$, this means that this locus is empty, that is, the only curve in $|\L|$ does not pass through any of the points of $S$.) Therefore
the condition that the nodes impose independent conditions on $|\L|$ can be rephrased as saying that the locus of curves in $|\L|$ passing through at least one node is an effective, simple normal crossing (SNC) divisor consisting of hyperplanes.}
\end{remark}

\renewcommand{\proofname}{Proof of Lemma \ref{lemma:r1}}

\begin{proof}
  By Serre duality we have that $H^0(\F^0_C) \iso Ext^1(\F^0_C,\omega_C)^{\vee}$. Applying 
 $\Shom_{\O_C}(-,\omega_C)$ to the bottom exact sequence of \eqref{eq:nod0r}, we find that $\Shext^1_{\O_C}(\F^0_C,\omega_C)=0$, so that
\[  Ext^1(\F^0_C,\omega_C)  \iso H^1(\Shom_{\O_C}(\F^0_C,\omega_C)) \]
by the local to global spectral sequence for $\Ext$. Moreover, applying $\Shom_{\O_C}(-,\omega_C)$ to \eqref{eq:r50}, we obtain the short exact sequence
\begin{equation} \label{eq:quagliarella}
\xymatrix{
0  \ar[r] & \Shom_{\O_C}(\F^0_C,\omega_C) \ar[r] & \N_{C/X} \ar[r] & \O_{S \cap C} 
\ar[r] & 0.
}
\end{equation}
Since $h^1(\O_X)=0$, the restriction map 
$ \xymatrix{ H^0(\L) \ar[r] & H^0(\N_{C/X})}$ is surjective
by
\begin{equation} \label{eq:sl}
\xymatrix{
0 \ar[r] &   \O_X \ar[r]   &  \L \ar[r]  &  \N_{C/X} \ar[r] &  0. 
}
\end{equation}
Hence
\[ \coker\Big\{ H^0(\N_{C/X}) \rightarrow H^0(\O_{S \cap C}) \Big\} = 
\coker\Big\{ H^0(\L) \rightarrow H^0(\O_{S \cap C}) \Big\}. \]
Taking cohomology of \eqref{eq:quagliarella} we therefore obtain
\[ H^1(\Shom_{\O_C}(\F^0_C,\omega_C)) \iso \coker\Big\{ H^0(\L) \rightarrow H^0(\O_{S \cap C}) \Big\} \+ H^1(\N_{C/X})\]
It follows that
\begin{eqnarray*}
  H^0(\F^0_C) & \iso & Ext^1(\F^0_C,\omega_C)^{\vee} \iso H^1(\Shom_{\O_C}(\F^0_C,\omega_C))^{\vee} \\
 & \iso & \coker\Big\{ H^0(\L) \rightarrow H^0(\O_{S \cap C}) \Big\}^{\vee} \oplus 
H^1(\N_{C/X})^{\vee} \\
& \iso & \coker\Big\{ H^0(\L) \rightarrow H^0(\O_{S \cap C}) \Big\}^{\vee} \oplus 
H^0(\N_{X/Y}\* \O_C),
\end{eqnarray*}
by Serre duality and \eqref{eq:r110} for the last isomorphism.
Therefore, $\iota_C$ is an isomorphism for all $C \in |\L|$ if and only if 
\[\xymatrix{
H^0(\L) \ar[r] & H^0(\O_{S \cap C}) 
}
\]
is surjective for all $C \in |\L|$, which precisely means that $S$ imposes independent conditions on $|\L|$.  
\end{proof}

\renewcommand{\proofname}{Proof}

The next result will be 
central in the rest of the paper.

\begin{prop} \label{lemma:r2}
  Assume that $\omega_Y \iso \O_Y$ and $h^1(\O_X)=0$.

We have  $H^0(C,\N_{C/X}) \iso H^0(C,\N_{C/Y})$ for all $C \in |\L|$ if and only if  
$S$ imposes independent conditions on $|\L|$ and the map $\gamma_C$ in \eqref{eq:r14} is an isomorphism for all $C \in |\L|$. 
\end{prop}

\begin{proof}
From the sequence \eqref{eq:diaF} we have that $H^0(\N_{C/X}) \iso H^0(\N_{C/Y})$ if and only if the connecting  homomorphism $\xymatrix{H^0(\F^0_C) \ar[r] &  H^1(\N_{C/X})}$ is injective. Since the spaces $H^0(C,\N_{X/Y}\* \O_C)$ and   $H^1(C,\N_{C/X})$ are  dual by 
\eqref{eq:r110}, this happens if and only if both the maps $\iota_C$ and $\gamma_C$ are isomorphisms. The result then follows from Lemma \ref{lemma:r1}.
\end{proof}

\section{Proof of Theorem \ref{holger}}  
\label{sec:mov}

In this section we prove Theorem \ref{holger}. The result is a consequence of \cite[Thms. 3.3 and 3.5]{ck} (and their proofs)
and \cite[Prop. 1.4 and Thm. 1.5]{kl1}, following the steps in the proof of \cite[Thm. 1.1]{kl2} almost ad verbatim except for an intermediate step where we apply 
Proposition \ref{lemma:r2} (cf. Remark \ref{rem:r5} below) and exploit the fact that conditions  (A3), (A4) and (A7) are stronger than the assumptions in
\cite[Thm. 1.1]{kl2}. Since \cite{kl2} has never been published, we give the whole proof of Theorem \ref{holger}.

The setting and assumptions are as in the introduction. We observe the following:

\begin{remark} \label{rem:r5}
  {\rm  By Proposition \ref{lemma:r2}, condition (A5) is equivalent to the condition}
    \begin{itemize}
    \item[(A5)'] {\rm  The set $S=\{\xi_1, \ldots, \xi_{\ell}\}$ imposes
 independent conditions on $|\L|$ and the map $\gamma_C$ in \eqref{eq:r14} is an isomorphism for all $C \in |\L|$.}
\end{itemize}
\end{remark}

\renewcommand{\proofname}{Proof of Theorem \ref{holger}}

\begin{proof}
Since $h^1(\O_X)=0$, we have that $|\L| \iso \PP^{\ell}$ is a connected component of $\Hilb X$. By (A5) we have that $|\L|$ is also a connected
component of $\Hilb Y$, and by (A6)  it has a
smooth neighborhood $\H \subset \Hilb P$.

Let
\[
\xymatrix{
         & \C_0   \ar[d]_{p_0} \ar[r]^{q_0}     & Y \\
& |\L| & 
}
\]
and
\[
\xymatrix{
         & \C   \ar[d]_{p} \ar[r]^{q}     & P \\
& \H & 
}
\]
be the universal curves. Let 
$\i$ be the ideal sheaf of $\C_0$ in $ |\L| \times Y$ and 
$\j$ be the ideal sheaf of $\C_0$ in $ |\L| \times P$.

Applying the functor
\[
F:={p_0}_*\circ \Shom_{\C_0} (-,\O_{\C_0})
\]
to the exact sequence
\[  
\xymatrix{
0 \ar[r] &  {q_0}^* \E ^{\vee} \ar[r] & \j/ \j^2    \ar[r] &  \i/ \i^2 \ar[r] & 0 
}
\]
of conormal sheaves and using the infinitesimal properties of Hilbert
schemes gives the exact sequence
\[ 
\xymatrix{
0 \ar[r] &  \T_{|\L|} \ar[r] & \T_{\H}\* \O_{|\L|} \ar[r] & {p_0}_*{q_0}^*\E 
\ar[r] & R^1 F (\i/ \i^2) \ar[r] & 0
}
\]
of $\O_{|\L|}$-modules, as $R^1 F (\j/ \j^2)=0$ by (A6).

Setting
\[
\Q:=R^1 F (\i\big/\i^2)
\]
we shorten the above to

\begin{equation}\label{eq:excess}
\xymatrix{
0 \ar[r] & \N_{|\L|/\H}  \ar[r] & {p_0}_*{q_0}^*\E  \ar[r]^{\hspace{0,4cm}\rho} & \Q \ar[r] & 0.
}
\end{equation}

By Remarks \ref{rem:snc} and \ref{rem:r5}, condition (A5) implies that the locus of curves in $|\L| \iso \PP^{\dim|\L|}$ passing through the node $\xi_i$ is a hyperplane
$\D_i \subset |\L|$, and furthermore, that the locus of curves in $|\L|$ passing through 
at least one node, satisfies
\begin{equation}
  \label{eq:snc}
  \D:= \cup_{i=1}^{\ell} \D_i \; \mbox{is simple normal crossing (SNC).}
\end{equation}

In fact, what we have done so far, is to go through the first part of the proof of \cite[Thm. 3.3]{ck}. Since the notation in \cite{ck} is rather involved, we 
include the following translation between our notation and the notation in 
\cite{ck}:

\vspace{0,5cm}

\begin{center}
\begin{footnotesize}
\begin{tabular}{|c|c|} \hline
Our notation & Notation in \cite{ck} \\ \hline \hline
$P$ & $P$        \\ \hline 
$Y=Y_0$  & $X_0$   \\ \hline 
$|\L|$ & $S'=S'(C)=I'$   \\ \hline 
$\H$   & $J'$   \\ \hline
$\C$ & $J$    \\ \hline
$\C_0$   & $S=I$   \\ \hline 
$p$, $q$ & $p$, $q$    \\ \hline 
${p_0}$, ${q_0}$ & $p_S$, $q_S$   \\ \hline 
$S$  & $\Xi$   \\ \hline 
$\ell$ & $r(S)$    \\ \hline
$\xi_i$ & $x^i_S$   \\ \hline 
$\D_i$, $\D$ & $D^i_{S'}$, $D_{S'}$    \\ \hline 
$X$ & $Y_S$   \\ \hline 
$\Q$ & $\Q_{S'}$    \\ \hline 
$p_*q^*\E$ & $\mathcal{V}$    \\ \hline 
\end{tabular}
\end{footnotesize}
\end{center}

\vspace{0,5cm}

In particular, the conditions in \cite[Thm. 3.3]{ck} are satisfied. Moreover, by (A4), the conditions in \cite[Thm. 3.5]{ck} are also satisfied, so that

\begin{equation} \label{eq:kons2}
 \Q  \iso \Omega^1_{|\L|}[\log \D],  
\end{equation}
the locally free sheaf of differentials with logarithmic poles along an SNC divisor consisting of $\ell$ hyperplanes (see e.g. \cite[\S 2]{ev} for the definition).

By \cite[Thm. 1.5]{kl1} and (A6) we have that $p_* q^* \E$ is locally free on $\H$ and $|\L|$ is the zero scheme of $p_* q^* s_0$, so that \eqref{eq:excess} identifies $\Q$ as the excess normal bundle to $p_* q^* s_0$, cf. \cite[\S 1.2]{kl1}.
Still by \cite[Thm. 1.5]{kl1}, the Hilbert scheme of the threefold
$Y_t := Z(s_0+ts)$, satisfies
\[ \Hilb Y_t \cap \H = Z(p_* q^* (s_0+ts)). \]
This already finishes the proof in the case $\dim |\L|=0$, by (A3). In the remaining cases it suffices to prove that $\rho({p_0}_*{q_0}^*s)$ 
will vanish at precisely ${\ell -2 \choose \dim|\L|}$ distinct 
points of $|\L|$, all corresponding to smooth, irreducible curves in $|\L|$. Indeed, 
by \cite[Prop. 1.4]{kl1} (and its proof), if $\rho({p_0}_*{q_0}^*s)$ has a 
{\it reduced and isolated} zero at a point $z \in |\L|$, then the scheme
$Z(p_* q^* (s_0+ts))$ has a reduced and isolated zero in any small enough complex analytic neighborhood of $z$ in $\H$ for $t >0$ sufficiently small. Since the curve corresponding to $z$ is smooth, the same holds true in this neighborhood of $z$. This will finish the proof of the theorem.

Assume then that $\dim |\L| >0$. The fact that $\ell \geq \dim |\L|+2$ in condition (A2) implies that the 
locally free sheaf
$\Q  \iso \Omega^1_{\PP^{\dim|\L|}}[\log \ell \; \PP^{\dim|\L|-1}]$ is globally 
generated (see e.g. \cite[Thm. 3.5]{dk} for a proof of this fact) and that
\begin{equation}
  \label{eq:ctop}
  \int_{|\L|} c_{top} (\Q) = {\ell -2 \choose \dim|\L|} >0. 
\end{equation}

Consider the standard exact sequence
\[ 
\xymatrix{
0  \ar[r] &  \Omega^1_{|\L|} \ar[r] & \Q \ar[r]^{\hspace{-0,5cm}\varepsilon=(\varepsilon_i)}
 &  \oplus_{i=1}^{\ell} \O_{\D_i} \ar[r] & 0
}
\]
into which $\Q  \iso \Omega^1_{|\L|}[\log \D]$ sits 
(see e.g. \cite[2.3(a)]{ev} or \cite[Prop. 2.3]{dk}), which is the same sequence as the one in \cite[Thm. 3.3]{ck}. By (A2) we can choose local coordinates in an analytic neighborhood of $0=\xi_i$
in $P$ such that 
\[s_0(x)= x_1 \+ \cdots \+ x_{r-4} \+ 
\Big(x_{r-3}^2 + \cdots + x_r^2\Big). \]
 Let $s \in H^0(P,\E)=\+_{i=1}^{r-3} H^0(P,\M_i)$ and write 
\[s(x)= f_1(x) \+ \cdots \+ f_{r-3}(x)\]
 in the same coordinates. In the proof of \cite[Thm. 3.3]{ck} it is shown that the composition
\[ 
\xymatrix{
H^0(P,\E)  \ar[r] &  H^0(|\L|, {p_0}_* {q_0}^*\E) \ar[r]^{\hspace{0,3cm}H^0(\rho)}  &  
H^0(|\L|, \Q) \ar[r]^{\hspace{-0,2cm}H^0(\varepsilon_i)} & H^0(\D_i, \O_{\D_i}) 
}
\]
is given by
\[ 
\xymatrix{
s \ar@{|->}[r] & f_{r-3}(0). 
}
\]
Therefore, the image of the composition
\[ 
\xymatrix{
H^0(P,\E)  \ar[r] &  H^0(|\L|, {p_0}_* {q_0}^*\E) \ar[r]^{\hspace{0,3cm}H^0(\rho)}  &  
H^0(|\L|, \Q) \ar[r]^{\hspace{-0,2cm}H^0(\varepsilon)} & H^0(\D, \+_{i=1}^{\ell}\O_{\D_i}) 
}
\]
equals the image of the natural restriction map 
\[ 
\xymatrix{H^0(P, \M_{r-3}) \ar[r] &  H^0(S, \M_{r-3} \* \O_S),}
\]
 which  
has codimension one by (A7). Since $h^0(\Q)=\ell-1$, as is well known (see e.g. \cite[Prop. 2.5]{dk}), and $H^0(|\L|,\Omega^1_{|\L|})=0$, the map $H^0(\varepsilon)$ is injective with image of codimension one. It follows that the composition morphism
\begin{equation}
  \label{eq:compo}
 \xymatrix{
H^0(P,\E)  \ar[r] &  H^0(|\L|, {p_0}_* {q_0}^*\E) \ar[r]^{\hspace{0,3cm}H^0(\rho)}  &  H^0(|\L|, \Q)  
} 
\end{equation}
is surjective. Thus, as $\Q$ is globally generated, it follows by \eqref{eq:ctop} that for the general section $s \in H^0(P, \E)$, we have that $\rho({p_0}_*{q_0}^*s)$ 
vanishes at precisely ${\ell -2 \choose \dim|\L|}$ {\it distinct} 
points of $|\L|$, all corresponding to smooth, irreducible curves in $|\L|$, by (A3), as desired. This finishes the proof of the theorem.
\end{proof}

\renewcommand{\proofname}{Proof}

\begin{remark}
  {\rm Looking more closely at the proof, we see that the condition that $\E$ splits as a direct sum of line bundles is not necessary: it would suffice (when $r \geq 6$) that
$\E = \F \+ \M$,
where $\F$ is a vector bundle of rank $r-4$ and $\M$ is a line bundle on $P$.
Writing $s_0 \in H^0(P,\E)$ as $s_0 = s_{0,\F} \+ s_{0,\M}$, with $s_{0,\F} \in
H^0(P,\F)$ and $s_{0,\M} \in H^0(P,\M)$, we would then have $Z:=Z(s_{0,\F})$ in the setting of Theorem \ref{holger}.}
\end{remark}

\section{$K3$ surfaces embedded in nodal Calabi-Yau complete intersection threefolds}

\label{sec:set}

The rest of the paper is devoted to proving Theorem \ref{result}.

 We first recall the well-known construction used
in \cite{cl1}, \cite{katz}, \cite{og} and \cite{EJS}
to embed a $K3$ surface into a nodal Calabi-Yau complete intersection 
($CICY$) threefold.

It is well known, and easily seen by adjunction,  that there are three types 
of $K3$ complete
intersection surfaces in projective space, namely the intersection types $(4)$
in $\PP^3$, $(2,3)$ in $\PP^4$ and $(2,2,2)$ in $\PP^5$.
Similarly, there are five types of $CICY$ threefolds in projective space, namely the intersection types $(5)$
in $\PP^4$, $(3,3)$ and $(4,2)$ in $\PP^5$, $(3,2,2)$ in $\PP^6$ and
$(2,2,2,2)$ in $\PP^7$.

Let $X$ be a $K3$ surface of degree $2\mu-2$ in $\PP ^{\mu}$ 
that is a complete intersection of type
$(a_1, \ldots ,a_{r-2})$ in some $\PP^r$, for $r \geq \mu$.
 We will always assume that 
\begin{equation}
  \label{eq:assa}
  a_i \geq 2 \; \; \mbox{for} \; \;  i \leq r-4 \; \; 
\mbox{and} \; \; a_{r-3} \geq a_{r-2},   
\end{equation}
 but we may have $a_{r-3}=1$ or $a_{r-2}=1$.

Let
\[ b_i = a_i \; \; \mbox{for} \; \; i=1, \ldots, r-2, \; \; 
\mbox{and} \; \; b_{r-3}=a_{r-3}+a_{r-2}. \]
Then each $b_i \geq 2$ and 
we can construct a Calabi-Yau threefold $Y$ that is a complete intersection 
of type
$(b_1, \ldots, b_{r-3})$ in $\PP^r$ as follows: Choose
generators $g_i$ of degrees $a_i$ for the ideal of $X$. 
So $X = Z(g_1, \ldots, g_{r-2})$. For general 
$\alpha_{ij} \in H^0(\PP^r,\O_{\PP^r}(b_i-a_j))$ define
\[ f_i:= \sum \alpha_{ij}g_j \]
and
\[ Y := Z(f_1, \hdots, f_{r-3}) \]
(here we follow \cite[Section 3]{kl1}, except for arranging indices in a different way). If the coefficient forms 
$\alpha_{ij}$ are chosen 
in a
sufficiently general way, $Y$ has only $\ell= (2\mu-2)a_{r-3}a_{r-2}$ ordinary double points  and they all lie on $X$.  
This can be checked using Bertini's theorem. In fact, the $\ell$ nodes are the intersection points of two general elements of $|\O_X(a_{r-3})|$ and $|\O_X(a_{r-2})|$ (distinct, when $a_{r-3}=a_{r-2}$).  As above, we denote the set of nodes by $S$.

Moreover, for general $\alpha_{ij}$, Bertini's theorem yields that the fourfold
\[ Z := Z(f_1, \hdots, f_{r-4}) \]
is smooth. (Note that $Z={\PP^r}$ if $r=4$.)

We  are therefore in the setting of Theorem \ref{holger} given in the introduction with $P=\PP^r$, 
\[ \E :=  \+ _{i=1}^{r-3} \O_{\PP^r}(b_i)  \]
and $\M_{r-3}:= \O_{\PP^r}(b_{r-3})= \O_{\PP^r}(a_{r-3}+a_{r-2})$. 
By construction, condition (A1) and the first part of condition (A2) are satisfied.

We refer to Table \ref{tabella}  for all values of
$a_j$, $b_i$, $\ell$, $\mu$ and $r$. (This is the same table as \cite[Table p.~201]{kl1}, except for one 
case, namely $(b_i)=(3,3)$, $(a_j)=(2,2,2)$, present in \cite[Table p.~201]{kl1} but absent in our table, because in this case none of the two cubic hypersurfaces will be smooth along $X$.)  

\begin{footnotesize}
\begin{table} \caption{Construction of $CICY$s} \label{tabella}
\begin{tabular}{|c|c|c|c|c|c|c|c|} \hline
$(b_i)$ & $(a_j)$ & $\mu$ & $r$ & $\ell$ & $\Sing Y$ & $a_{r-3}$ & $a_{r-2}$ 
\\ \hline \hline
$(5)$  & $(4,1)$  &  $3$ & $4 $ & $16$ & $X \cap Z(\alpha_{11}, \alpha_{12})$ & $4$ & $1$ 
\\ \hline
$(5)$  & $(3,2)$  &  $4$ & $4$ & $36$ & $X \cap Z(\alpha_{11}, \alpha_{12})$ & $3$ & $2$ 
\\ \hline
$(4,2)$  & $(4,1,1)$  &  $3$ & $5$ &  $4$ & 
$X \cap Z(\alpha_{22}, \alpha_{23})$ & $1$ & $1$ 
\\ \hline
$(2,4)$  & $(2,3,1)$  &  $4$ & $5$ & $18$ & 
$X \cap Z(\alpha_{11}, \alpha_{12}\alpha_{23}-\alpha_{13}\alpha_{22})$ & $3$ & $1$ 
\\ \hline
$(2,4)$  & $(2,2,2)$  &  $5$ & $5$ & $32$ & 
$X \cap Z(\alpha_{21}\alpha_{12}-\alpha_{22}\alpha_{11}, \alpha_{21}\alpha_{13}-\alpha_{23}\alpha_{11})$ & $2$ & $2$ 
\\ \hline
$(3,3)$  & $(3,2,1)$  &  $4$ & $5$ & $12$ & 
$X \cap Z(\alpha_{21}\alpha_{12}-\alpha_{22}\alpha_{11}, \alpha_{21}\alpha_{13}-\alpha_{23}\alpha_{11})$ & $2$ & $1$ 
\\ \hline
$(3,2,2)$  & $(3,2,1,1)$  &  $4$ & $6$ & $6$ & 
$X \cap Z(\alpha_{22}\alpha_{33}-\alpha_{23}\alpha_{32}, 
\alpha_{22}\alpha_{34}-\alpha_{24}\alpha_{32})$ & $1$ & $1$ 
\\ \hline
$(2,2,3)$  & $(2,2,2,1)$  &  $5$ & $6$ & $16$ & $X \cap Z(linear, quadratic)$ & $2$ & $1$ 
\\ \hline
$(2,2,2,2)$  & $(2,2,2,1,1)$  &  $5$ & $7$ & $8$ & $X \cap Z(linear, quadratic)$ & $1$ & $1$ 
\\ \hline
\end{tabular}
\end{table}
\end{footnotesize}

Assume now that $X$ carries a line bundle $\L$
 such that  the general element of $|\L|$ is a smooth, 
irreducible curve of degree $d$ and genus $g$. It is well known, and easily seen, that such a line bundle satisfies
\begin{equation}
  \label{eq:cohl}
  \L^2=2g-2, \; h^0(\L)=g+1, \; h^1(\L)=h^2(\L)=0,
\end{equation}
see e.g. \cite{S-D}.

In the next section we will finish the proof of Theorem \ref{result} by applying Theorem \ref{holger} and a result guaranteeing the existence of the line bundle $\L$, cf. Theorem \ref{thm:bncurves}. In particular, we will verify that conditions (A2)-(A7) are satisfied (under certain numerical conditions, giving the different constraints in Theorem \ref{result}). In the rest of this section we will give some results that will be needed for the verifications of the conditions (A4), (A5) and 
(A7), where we in the case of (A5) will consider the equivalent condition 
(A5)' from Remark \ref{rem:r5}.

 The following result, which is folklore among experts on $K3$ surfaces, will be needed to verify that condition (A4) is satisfied.

\begin{lemma} \label{lemma:curvesonk3}
  Let $\L$ be a line bundle on a $K3$ surface $X$ such that $|\L|$ contains a smooth, irreducible curve. Assume that $x$ is a point on $X$ satisfying
  \begin{itemize}
  \item[(i)] $x$ is not contained in any smooth rational curve $\Gamma$ on $X$ satisfying $\Gamma.\L=0$;
  \item[(ii)] if $\L^2=0$, then $x$ is not a singular point of any 
fiber of the elliptic pencil $|\L|$.
\end{itemize}

Then, the general element in $|\L \* \I_x|$ (if nonempty) is nonsingular at $x$.
\end{lemma}

\begin{proof} We can assume that $g=\dim |\L|>0$. 
We then have that $|\L|$ is base point free, 
see \cite[Thm. 3.1]{S-D}, and $\dim |\L|=g$ and $\L^2=2g-2 \geq 0$ by 
\eqref{eq:cohl}.

If $\L^2=0$ it is well known that $|\L|$ is an elliptic pencil, see 
\cite[Prop. 2.6]{S-D}, so that we are done by (ii). (Note that $|\L \* \I_x|$ has only one element.)

If $\L^2 >0$ we consider the morphism $\varphi_{\L}:X \khpil \PP^g$ and its Stein factorization
\[
 \xymatrix{
X  \ar[r]^{\hspace{-0,1cm}\alpha} &  X' \ar[r]^{\hspace{-0,2cm}\beta}  & 
\varphi_{\L}(X) 
} 
\]
Then $\alpha$ is the contraction of the finitely many smooth, rational curves $\Gamma$ satisfying $\Gamma.\L=0$ \cite[(4.2)]{S-D} and $\beta$ is finite of degree one or two \cite[(4.1)]{S-D}. It follows from hypothesis (i) that the base scheme of $|\L \* \I_x|$ is finite of length at most two, whence is curvilinear. Therefore 
the general member of $|\L \* \I_x|$ is smooth and irreducible by Bertini. 
\end{proof}

By this lemma, condition (A4) is satisfied for general $\alpha_{ij}$. Indeed, the reduced and irreducible curves $\Gamma$ on $X$ satisfying $\Gamma.\L=0$ are only finitely many by standard results, and the singular points of fibers in an elliptic pencil are also finitely many. Therefore, since the nodes are the points of intersection of two general (distinct) elements of $|\O_X(a_{r-3})|$ and $|\O_X(a_{r-2})|$, we can make sure that the nodes satisfy the conditions (i) and (ii) in the lemma.

The next result shows that condition (A7) is satisfied.
 
\begin{lemma} \label{lemma:a7}
 For general $\alpha_{ij}$, the image of the natural restriction map 
\[  
\xymatrix{
H^0({\PP^r}, \O_{{\PP^r}}(a_{r-3}+a_{r-2})) \ar[r] & H^0(S, \O_S(a_{r-3}+a_{r-2})) \iso \CC^{\ell} 
}
\]
has codimension one. 
\end{lemma}

\begin{proof}
Recall that $S = H_1 \cap H_2$, for general members
$H_1 \in |\O_X(a_{r-3})|$ and $H_2 \in |\O_X(a_{r-2})|$.
We note that the map above factorizes through the natural restrictions to $X$ and $H_1$. The restriction map 
$H^0({\PP^r}, \O_{{\PP^r}}(a_{r-3}+a_{r-2})) \khpil H^0(X, \O_{X}(a_{r-3}+a_{r-2}))$ is surjective because $X$ is a complete intersection. The restriction map 
$H^0(X, \O_{X}(a_{r-3}+a_{r-2})) \khpil H^0(H_1, \O_{H_1}(a_{r-3}+a_{r-2}))$ is surjective from the cohomology of 
\[  
\xymatrix{ 0 \ar[r] &   \O_X(a_{r-2}) \ar[r] &  \O_X(a_{r-3}+a_{r-2}) \ar[r] &   \O_{H_1}(a_{r-3}+a_{r-2}) \ar[r] &   0, 
}
\]
as $h^1(\O_X(a_{r-2}))=0$. Finally, the cokernel of the map 
$H^0(H_1, \O_{H_1}(a_{r-3}+a_{r-2})) \khpil H^0(S, \O_{S}(a_{r-3}+a_{r-2}))$
is one-dimensional from the cohomology of 
\[ \xymatrix{ 0 \ar[r] &   \O_{H_1}(a_{r-3}) \ar[r] &  \O_{H_1}(a_{r-3}+a_{r-2}) \ar[r] &   \O_{S}(a_{r-3}+a_{r-2}) \ar[r] &   0,
}
 \]
as $\O_{H_1}(a_{r-3}) \iso \omega_{H_1}$. This concludes the proof.
\end{proof}

The following result gives criteria for the first of the two conditions in (A5)'  to hold:

\begin{lemma} \label{lemma:tappato}
 Assume that the $\alpha_{ij}$ are general, that $h^0(X,\L \* \O_X(-a_{r-2}))=0$ and that
\[
a_{r-2}(2a_{r-3}-a_{r-2})(\mu-1) \geq
\begin{cases} 
g+2 & \; \mbox{if $a_{r-3} \neq a_{r-2}$;} \\
g+1 & \; \mbox{if $a_{r-3}=a_{r-2}$.}
\end{cases} 
\]

Then $S$ imposes independent conditions on $|\L|$.
\end{lemma}

\begin{proof}
Fix any smooth, irreducible 
$H_0 \in |\O_X(a_{r-2})|$.  
Consider the incidence scheme
\[ 
W:= \Big\{ (\eta',\eta) \; | \; \eta' \subset \eta \Big \} \subset \Sym^{g+1}(H_0) \x |\O_{H_0}(a_{r-3})|
\]
(note that $\deg \O_{H_0}(a_{r-3})=\ell=2a_{r-3}a_{r-2}(\mu-1) \geq g+1$ by our conditions). Let
\[ \xymatrix{  \pi_1: W \ar[r] &  \Sym^{g+1}(H_0)} \; \; \mbox{and} \; \; \xymatrix{ \pi_2: W \ar[r] &  |\O_{H_0}(a_{r-3})|} \]
denote the projections. One easily computes that 
\[
\dim |\O_{H_0}(a_{r-3})|= 
\begin{cases} 
a_{r-2}(2a_{r-3}-a_{r-2})(\mu-1) -1& \; \mbox{if $a_{r-3} \neq a_{r-2}$;} \\
a_{r-2}(2a_{r-3}-a_{r-2})(\mu-1) & \; \mbox{if $a_{r-3}=a_{r-2}$,}
\end{cases} 
\]
so that $\dim |\O_{H_0}(a_{r-3})| \geq g+1$ by our assumptions. It follows that 
$\pi_1$ is surjective. Since, for any $\eta' \in \Sym^{g+1}(H_0)$, we have
$\pi_1^{-1}(\eta') \iso |\O_{H_0}(a_{r-3}) \* \I_{\eta'/H_0}|$, it follows that $W$ is irreducible. 

For general $\eta' \in \Sym^{g+1}(H_0)$, we have 
\begin{equation}
  \label{eq:vuoto}
 |\L \* \I_{\eta'/X}| =\emptyset. 
\end{equation}
Indeed, for general $x_1 \in H_0$, we have $\dim|\L \* \I_{x_1/X}|=\dim|\L|-1$, as $H_0$ cannot be a base component of $|\L|$. Proceeding inductively, having picked general distinct $x_1, \ldots, x_{i} \in H_0$ for some $i \in \{1, \ldots, g\}$ with
 $\dim|\L \* \I_{x_1/X}\* \cdots \* \I_{x_i/X}|=\dim|\L|-i \geq 1$, as $\dim|\L|=g$, we have that the base locus of $|\L \* \I_{x_1/X}\* \cdots \* \I_{x_i/X}|$ does not contain the whole of $H_0$, as $h^0(\L-H_0)=h^0(\L \* \O_X(-a_{r-2}))=0$ by our assumption. Therefore, for any $x_{i+1} \in H_0$ outside of this base locus, we have that $\dim|\L \* \I_{x_1/X}\* \cdots \* \I_{x_{i+1}/X}|=\dim|\L|-(i+1)$ and we can continue, proving \eqref{eq:vuoto}. 

Let now
\[ 
W_{\L}:= \Big\{ (\eta',\eta) \in W \; \Big| \; |\L \* \I_{\eta'/X}| \neq \emptyset \Big\} \sub W.
\]
Since $\pi_1$ is surjective, we have $W_{\L} \subsetneqq W$ by \eqref{eq:vuoto}. Since $\pi_2$ is finite, we get that $\pi_2(W_{\L}) \subsetneqq |\O_{H_0}(a_{r-3})|$. Since $h^1(\O_X(a_{r-3}-a_{r-2}))=0$, we conclude from the short exact sequence 
\[  
\xymatrix{
0 \ar[r] &  \O_X(a_{r-3}-a_{r-2}) \ar[r] & \O_X(a_{r-3}) \ar[r] & \O_{H_0}(a_{r-3}) 
\ar[r] &  0,
}
\]
that for the general member $H_1 \in |\O_X(a_{r-3})|$, the $\ell$ distinct points
$H_1 \cap H_0$ have the property that no $g+1$ of them lie on any member of $|\L|$. It follows that $S$ imposes independent conditions on $\L$, cf. Remark \ref{rem:snc}.
\end{proof}

Finally, we will consider the second of the two conditions in (A5)'.

We first recall the following result:

\begin{lemma}{\rm(\cite[Lemma 1.10]{kl1})} \label{lemma:kl110}
  Let $X \subset \PP^r$ be a smooth projective $K3$ surface with a line bundle $\L$ such that $\L$ and $\O_X(1)$ are independent in $\Pic X$ and 
$|\L|$ contains a smooth, irreducible curve. Then, for all $C \in |\L|$, the composition
\[  
\xymatrix{
\varphi: H^0(X, \N_{X/\PP^{r}}) \ar[r] & H^0(X, \N_{X/\PP^r} \* \O_C) \ar[r] & H^1(C, \N_{C/X})
}
\]
of the restriction with the connecting homomorphism arising from 
\eqref{eq:stand} with $P=\PP^r$ is surjective. Furthermore, $\ker \varphi$ is independent of $C \in |\L|$.
\end{lemma}

 We now prove that, for 
general $\alpha_{ij}$, the second of the two conditions in (A5)' from Remark \ref{rem:r5} holds if $\L$ and $\O_X(1)$ are independent in 
$\Pic X$. 
The proof is taken basically ad verbatim from the last part of the proof of \cite[Thm. 3.5]{kl1} and we give it not only for the sake of completeness, but also because our statement is different and the proof of \cite[Thm. 3.5]{kl1} has a gap (cf. Remark \ref{rem:hf} below).

\begin{prop} \label{thm:gap}
Suppose that the $\alpha_{ij}$ are general and that the line bundles $\L$ and $\O_X(1)$ are independent in 
$\Pic X$.  Then the map $\gamma_C$
in \eqref{eq:r14} is an isomorphism for all $C \in |\L|$. 
\end{prop} 

\renewcommand{\proofname}{Proof (following the proof of \cite[Thm. 3.5]{kl1})}

\begin{proof}
We consider the commutative diagram \eqref{eq:nod0r} with $P=\PP^r$.
By definition of $\gamma_C$ and the fact that $h^1(\N_{C/X})=1$ for any $C \in |\L|$ by \eqref{eq:sl} and \eqref{eq:cohl}, we just need to show that the connecting homomorphism of the left hand vertical exact sequence in \eqref{eq:nod0r} is nonzero.
 A diagram chase in \eqref{eq:nod0r} reduces this to proving that for general $\alpha_{ij}$ and all $C \in |\L|$, we have
\begin{equation} \label{eq:claim}
 \delta_C(\ker H^0(\Phi_C)) \neq 0,
\end{equation}
where 
\[  
\xymatrix{
\delta_C: H^0(C, \N_{X/\PP^r} \* \O_C) \ar[r] & H^1(C,\N_{C/X})
}
\] 
is the connecting homomorphism of the right hand vertical exact sequence in \eqref{eq:nod0r}
and $\Phi_C$ is given in \eqref{eq:nod0r}.

As in \cite{kl1}, we define
$\A := \+ _{j=1}^{r-2} \O_{\PP^r}(a_j)$ and denote by $\A_{X}$ its restriction to $X$.

Now we have natural isomorphisms 
\begin{equation}
  \label{eq:kl34}
  \N_{X/\PP^r} \iso \A_X = \+ _{j=1}^{r-2} \O_{X}(a_j) \; \; \mbox{and} \; \; 
\N_{Y/\PP^r} \* \O_X \iso \E  \* \O_X \iso \+ _{i=1}^{r-2} \O_{X}(b_i), 
\end{equation}
under which the map
\[  
\xymatrix{
\+ _{j=1}^{r-2} H^0(\O_{C}(a_j)) \iso H^0(\N_{X/\PP^r} \* \O_C) \ar[r]^{H^0(\Phi_C)} & H^0(\N_{Y/\PP^r} \* \O_C) \iso \+ _{i=1}^{r-2} H^0(\O_{C}(b_i))
}
\]
is given by the matrix $(\overline{\alpha}_{ij})$, where 
$\overline{\alpha}_{ij}$ is the restriction of $\alpha_{ij}$ to $C$ (cf. \cite[Lemma 3.4]{kl1}, where $\B :=  \+ _{i=1}^{r-3} \O_{{\PP^r}}(b_i) = \E$). 

To prove \eqref{eq:claim}, use Gaussian elimination to find a generator $N=N(\alpha_{ij})$ of the null-space of the linear map 
$\xymatrix{ \CC^{r-2} \ar[r]^{(\alpha_{ij})} & \CC^{r-3}}.$ Keeping track of degrees, this can be done in such a way that the $i$th coordinate of $N$ is of degree $a_i$. For example, if $(b_i)=(4,2)$ and $(a_j)=(3,2,1)$, the vector 
$N=(\alpha_{12}\alpha_{23}-\alpha_{13}\alpha_{22}, -\alpha_{11}\alpha_{23},
\alpha_{11}\alpha_{22})$. Since $N$ is well-defined only up to scalar, we view it as a line in $H^0(\PP^r, \A)$ or a point of $\PP(H^0(\PP^r, \A))$.
Now each term of each coordinate of $N$ is a term of a determinant of 
$(\alpha_{ij})$, so that as $(\alpha_{ij})$ varies, $N$ hits a multiple of each element of the form $(0, \ldots, 0, \lambda_1 \lambda_2 \cdots \lambda_{a_i},0, \ldots, 0)$, where the $\lambda_k$s are of degree $1$. Since the image of the Segre embedding of $\PP(H^0(\PP^r,\O(1)))^{\x a_i}$ in $\PP(H^0(\PP^r,\O(a_i)))$ is nondegenerate, the linear span of the image of $N$ includes the spaces
\[ 0 \+ \cdots \+ 0 \+ H^0(\PP^r,\O(a_i)) \+ 0 \+ \cdots \+ 0 \]
and hence all of $\+_i H^0(\PP^r, \O(a_i))$.

Now \eqref{eq:claim} follows from Lemma \ref{lemma:kl110} since 
$H^0(\PP^r, \A) \khpil H^0(X, \A_X)$ is surjective for any complete intersection $X$.   
\end{proof}

\renewcommand{\proofname}{Proof}

\begin{remark} \label{rem:hf}
{\rm In \cite[Thm. 3.5]{kl1} it is claimed that $H^0(\N_{C/X}) \iso H^0(\N_{C/Y})$ for all $C \in |\L|$  under the same hypotheses as 
in Proposition \ref{thm:gap}. 
However, there is a mistake in the proof: 
it is claimed that
there is a (not necessarily surjective) map $\mu: \N_{C/Y} \khpil \N_{X/Y} \* \O_C$ fitting into the left bottom corner of \eqref{eq:nod0r}. While this is true from \eqref{eq:diaF}
if $C \cap S = \emptyset$, because then
$\F^0_C = \N_{X/Y} \* \O_C$ by \eqref{eq:r110} and \eqref{eq:r50},
the same does not hold if $C$ passes through some nodes of $Y$. In this case, $\N_{X/Y} \* \O_C$ is a {\it proper} subsheaf of $\F^0_C$   and one cannot conclude that there is a map as $\mu$. 

Moreover, the gap in the proof of \cite[Thm. 3.5]{kl1} also influences the proof of \cite[Corollary 3.6]{kl1}, which in fact does not hold. Indeed, if 
$H^0(\N_{C/X}) \iso H^0(\N_{C/Y})$ for all $C \in |\L|$, then by 
\eqref{eq:diaF}, \eqref{eq:r50} and the fact that $h^1(\N_{C/X})=h^1(\O_C(C))=1$ for every $C \in |\L|$,  we must have
\[ h^1(\N_{C/Y})=h^1(\F^0_C)= g - \sharp(S \cap C). \]
It therefore follows that \cite[Corollary 3.6]{kl1} does not hold for the curves $C \in |\L|$ passing through at least one of the nodes, as it is claimed there that $H^1(C,\N_{C/Y}) \iso H^1(C, \O_C)$ for {\it any} $C \in |\L|$, and it is well known that $h^1(\O_C)=g$ for any $C \in |\L|$, cf. e.g. 
\cite[Lemma 1.9(1)]{kl1}. Unfortunately, \cite[Corollary 3.6]{kl1} is used in 
a crucial way in the proof of the main existence result \cite[Thm. 1]{kl1} (more precisely, it is used in the proof of \cite[Prop. 3.9]{kl1}).
}
\end{remark}

\begin{remark} \label{rem:hh}
{\rm It was noted by the referee of the original version of \cite{kn2} that the same gap as in the proof of \cite[Thm. 3.5]{kl1} also appears in \cite[Example 4.3]{ck}. It may be instructive to have a closer look at this as well and to see how Proposition \ref{lemma:r2} can be used to fill this gap and thus show that the conclusions in \cite[Example 4.3]{ck} are correct.

Clemens and Kley consider, with different notation, the two upper cases in 
Table \ref{tabella}, that is, either
\begin{equation}
  \label{eq:es1}
  g_1, \alpha_{12} \in H^0(\PP^4,\O_{\PP^4}(4)) \; \mbox{and} \; 
g_2, \alpha_{11} \in H^0(\PP^4,\O_{\PP^4}(1))
\end{equation}
or
\begin{equation}
  \label{eq:es2}
  g_1, \alpha_{12} \in H^0(\PP^4,\O_{\PP^4}(3)) \; \mbox{and} \; 
g_2, \alpha_{11} \in H^0(\PP^4,\O_{\PP^4}(2))
\end{equation}
(respectively cases (4.3.1) and (4.3.2) in their example), so that both the $K3$ surface
\[ X:=Z(g_1,g_2) \]
(which they call $Y$)
and the del Pezzo surface
\[ X':=Z(g_2, \alpha_{11}) \]
(which they call $S$) are smooth, and such that the quintic threefold
\[ Y:= Z(\alpha_{11}g_1+\alpha_{12}g_2) \]
(which they call $X_0$)
has only ordinary nodes ($16$ and $36$ respectively), given by
\[ S=Z(g_1,g_2,\alpha_{11},\alpha_{12}) \subset X \cap X'. \]
The authors claim that for any curve $C' \subset X'$, one has 
\begin{equation}
  \label{eq:wrong}
  H^0(\N_{C'/X'}) \iso H^0(\N_{C'/Y}).
\end{equation}
 To prove \eqref{eq:wrong} they again state the existence of an exact sequence 
\[  
\xymatrix{ 0 \hpil H^0(C', \N_{C'/X'}) \ar[r] &  H^0(C', \N_{C'/Y}) \hpil 
H^0(C',\N_{X'/Y} \* \O_{C'}).} \]
As mentioned above, one cannot a priori conclude the existence of such a sequence. In fact, as we will see now, \eqref{eq:wrong} does not hold in certain cases. 

By Proposition \ref{lemma:r2}  a necessary and sufficient condition for \eqref{eq:wrong} to hold for all curves $C'$ in the linear system $|C'|$ is that the nodes $S$  impose independent conditions on $|C'|$, since $X'$ is a del Pezzo surface,  so that $h^1(\N_{C/X'})=0$ for all $C \in |C'|$ and $\gamma_C$ is automatically surjective.

Now in case \eqref{eq:es1}, let $C'$ be the quartic $X' \cap Z(g_1)$ on $X'$.
Then $C'$ pass through all $16$ nodes but $|\O_{X'}(C')|=|\O_{\PP^2}(4)|$ is $14$-dimensional. Hence the nodes $S$ do not impose independent conditions on $|C'|$,   so that 
\eqref{eq:wrong} does not hold for all curves in $|C'|$
in this case. In fact, it does not hold 
when one considers linear systems $|\O_{\PP^2}(d)|$ with $d \geq 4$. This example also shows that the assertion ``If the $g_i$ and the $\alpha_{ij}$ are sufficiently general, the divisor of curves passing through at least one node is a simple-normal crossing divisor consisting of hyperplanes'' in
\cite[Example 4.3]{ck} does not hold if $d \geq 4$, cf. Remark \ref{rem:snc}.

We now show, however, that $S$ does impose independent conditions on $|\O_{\PP^2}(d)|$ when $d \leq 3$, using the same argument as in the proof of Lemma \ref{lemma:tappato}.

Remember that the $16$ nodes of $Y$ are the intersection points of two general elements of $|\O_{\PP^2}(4)|$.
Fix any smooth, irreducible $H_0 \in |\O_{\PP^2}(4)|$ 
and consider the incidence scheme
\[ 
W:= \Big\{ (\eta',\eta) \; | \; \eta' \subset \eta \Big \} \subset \Sym^{\frac{1}{2}d(d+3)+1}(H_0) \x |\O_{H_0}(4)|
\]
with projections $\pi_1$ and $\pi_2$ onto the first and second factor, respectively (recall that $\dim |\O_{\PP^2}(d)|=\frac{1}{2}d(d+3)$). Since 
$\dim |\O_{H_0}(4)|=14 \geq \frac{1}{2}d(d+3)+1$ when $d \leq 3$, we have that 
$\pi_1$ is surjective. Since, for any $\eta' \in \Sym^{\frac{1}{2}d(d+3)+1}(H_0)$, we have
$\pi_1^{-1}(\eta') \iso |\O_{H_0}(4) \* \I_{\eta'/H_0}|$, it follows that $W$ is 
irreducible. 

Now using the fact that $|C' \* \O_{\PP^2}(-4)|=|\O_{\PP^2}(d-4)| = \emptyset$ as $d \leq 3$, one can easily prove that $|\O_{\PP^2}(d) \* \I_{\eta'/X}| =\emptyset$
for general $\eta' \in \Sym^{\frac{1}{2}d(d+3)+1}(H_0)$, as in the proof of Lemma \ref{lemma:tappato}. Now $\pi_1$ is surjective, so that letting
\[ 
W_{\L}:= \Big\{ (\eta',\eta) \in W \; | \; \eta' \subset C' \; \mbox{for some} \; C' \in |\O_{\PP^2}(4)| \Big\} \sub W,
\]
this implies that $W_{\L} \subsetneqq W$. Since $\pi_2$ is finite, we get that $\pi_2(W_{\L}) \subsetneqq |\O_{H_0}(4)|$. Since $H^0(\O_{\PP^2}(4))$ surjects 
onto $H^0(\O_{H_0}(4))$, we see that for the general member $H_1 \in |\O_{\PP^2}(4)|$, the $16$ distinct points
$H_1 \cap H_0$ have the property that no $\frac{1}{2}d(d+3)+1$ of them lie on any member of $|\O_{H_0}(d)|$, for $d \leq 3$. Therefore, for general choices of $g_1$ and $g_2$, \eqref{eq:wrong} is 
fulfilled for all curves $C' \in |\O_{\PP^2}(d)|$ when $d \leq 3$

The cases where $d \leq 3$ are in fact the cases where 
$H^1(\N_{C'/{\mathbb P}^4})=0$ for all $C' \in |\O_{\PP^2}(d)|$, another condition needed to apply the results in 
\cite{ck}, so that at the end, these are precisely the cases in \cite[Example 4.3]{ck} where the authors compute the number of curves moving in a general deformation of the threefold. This shows that the conclusions in \cite[Example 4.3]{ck}
in the case \eqref{eq:es1} are in fact correct.  

In the same way one can show that the applications in the case \eqref{eq:es2}
are correct.
}
\end{remark}

\section{Proof of Theorem \ref{result}} 
\label{sec:con}

 In this section we apply Theorem \ref{holger} and the results from the  previous section to 
construct smooth, isolated curves in general $CICY$s of each intersection type,
thus proving Theorem \ref{result}.

Recall the following result from \cite{mori} and \cite{kn1}, which guarantees the existence of smooth curves of certain degrees and genera on the three types of complete intersection $K3$ surfaces:

\begin{theorem} \label{thm:bncurves}
  Let $d>0$ and $g \geq 0$ be integers. Then:
  \begin{itemize}
    \item[(i)] There exists a smooth quartic
    surface $X$ in
    $\PP^3$ containing a smooth, irreducible curve $C$ of degree $d$ and genus $g$ such that $\O_X (1)$ and
    $\O_X (C)$ are independent in $\Pic X$ if and
    only if $g<d^2/8$ and $(d,g)
    \not = (5,3)$.

    \item[(ii)] There exists a $K3$ surface $X$ of type $(2,3)$ in  $\PP^4$
    containing a smooth, irreducible curve $C$ of degree $d$ and genus $g$ 
such that $\O_X (1)$ and
    $\O_X (C)$ are independent in $\Pic X$ 
if and
    only if $g=d^2/12+1/4$ or $g<d^2/12$ and $(d,g)\not = (7,4)$. 

    \item[(iii)] There exists a $K3$ surface $X$ of type $(2,2,2)$ in  $\PP^5$
    containing a smooth, irreducible curve $C$ of degree $d$ and genus $g$  
such that $\O_X (1)$ and
    $\O_X (C)$ are independent in $\Pic X$ if and
    only if $g=d^2/16$
    and $d \eqv 4  \hs (\mod 8)$, or $g<d^2/16$ and $(d,g )
    \not = (9,5)$.
\end{itemize}
\end{theorem}

For $X$ and $C$ as in the theorem, we let $\L:=\O_X(C)$.

In the setting described in the previous section, we now want to find out under which circumstances the conditions (A1)-(A7) in Section \ref{sec:mov} are satisfied, so that we can apply Theorem \ref{holger}.

\begin{prop} \label{prop:checkcond}
  Under the contraints given by Theorem \ref{thm:bncurves}, assume that the 
  $\alpha_{ij}$ are general. If
\begin{equation}
  \label{eq:min}
  d \leq 2a_{r-2}(\mu-1) \; \; \mbox{or} \; \; da_{r-2} > a_{r-2}^2(\mu-1)+g 
\end{equation}
and 
\begin{equation}
  \label{eq:A3}
a_{r-2}(2a_{r-3}-a_{r-2})(\mu-1) \geq
\begin{cases} 
g+2 & \; \mbox{if $a_{r-3} \neq a_{r-2}$;} \\
g+1 & \; \mbox{if $a_{r-3}=a_{r-2}$},
\end{cases} 
\end{equation}
the conditions (A1)-(A7) are satisfied.
\end{prop}

\begin{proof}
Conditions (A1) and (A3) are obviously satisfied 
(as is the first part of (A2)). Condition (A4) is satisfied by Lemma \ref{lemma:curvesonk3} and the lines following the proof of that lemma. 
Condition (A7) is satisfied by Lemma \ref{lemma:a7}. 

The remaining three conditions (A2), (A5) and (A6) will be responsible for the 
numerical conditions \eqref{eq:min} and \eqref{eq:A3}.

Condition 
(A2) is satisfied whenever
\begin{equation}
  \label{eq:A2}
  \ell \geq g+2,
\end{equation}
since $\dim |\L|=g$ by \eqref{eq:cohl}, where $\ell= 2a_{r-3}a_{r-2}(\mu-1)$ is the number of nodes of $Y$. We see that condition \eqref{eq:A3}
implies condition \eqref{eq:A2}.

We next consider condition (A6). Lemma \ref{lemma:kl110} together with \eqref{eq:stand} and the fact that 
 $\N_{X/{\PP^r}} \* \O_C = \oplus \O_C(a_j)$ by \eqref{eq:kl34} yields 
\[
H^1(C,\N_{C/\PP^r}) \iso H^1(C, \N_{X/\PP^r} \* \O_C) \iso \+ H^1(C, \O_C(a_j)). 
\]
Noting that $a_{r-2}=\min\{a_j\}$ by Table \ref{tabella}, we get in particular that 
$H^1(C, \N_{C/\PP^r})=0$  if and only if $h^1(C, \O_C(a_{r-2}))=0$.
By \cite[Prop. 1.3]{kn1} this is achievable for all $C \in |\L|$ if and only if 
\[
  d \leq 2a_{r-2}(\mu-1) \; \; \mbox{or} \; \; da_{r-2} > a_{r-2}^2(\mu-1)+g, 
\]
which is condition \eqref{eq:min}.

Therefore, (A6)  is satisfied if and only if \eqref{eq:min} 
holds.

Finally we consider condition (A5), or
equivalently, by Remark \ref{rem:r5}, condition (A5)'.

By Proposition \ref{thm:gap} and the fact that
$\L$ and $\O_X(1)$ are independent in $\Pic X$ by Theorem \ref{thm:bncurves}, the second of the two conditions in (A5)' is satisfied. 

Next we note from the cohomology of
\[ 
\xymatrix{
0  \ar[r] &  \L^{\vee} \ar[r] &  \O_X \ar[r] & \O_C \ar[r] & 0
}
\]
twisted by $\O_X(a_{r-2})$, Kodaira vanishing and Serre duality, that
\[ h^0(X, \L \otimes \O_X(-a_{r-2})) = h^1(\O_C(a_{r-2})), \]
so that also $h^0(X, \L \otimes \O_X(-a_{r-2}))=0$ if condition \eqref{eq:min}
holds, as we have just seen. Thus, by Lemma \ref{lemma:tappato}, the first of the two conditions  in (A5)' is satisfied whenever \eqref{eq:min} holds together with the condition 
\[
a_{r-2}(2a_{r-3}-a_{r-2})(\mu-1) \geq
\begin{cases} 
g+2 & \; \mbox{if $a_{r-3} \neq a_{r-2}$;} \\
g+1 & \; \mbox{if $a_{r-3}=a_{r-2}$,}
\end{cases} 
\]
which is condition \eqref{eq:A3}.

To summarize, conditions (A1)-(A7) are satisfied  for general $\alpha_{ij}$ whenever \eqref{eq:min} and
 \eqref{eq:A3} hold, subject to the constraints given by Theorem \ref{thm:bncurves}.
\end{proof}

Now to obtain Theorem \ref{result} we just apply Theorem \ref{holger} taking into account the numerical conditions given in Theorem \ref{thm:bncurves} and Proposition \ref{prop:checkcond} 
for the various complete intersection types. Let us briefly explain how it works in case (a), that is the case of $Y$ a quintic threefold in $\PP^4$. We then have $r=4$ and $(b_i)=(b_1)=5$.

Looking at Table \ref{tabella}, there are two choices for the $K3$ surface $X$, and the values of $a_j, \mu, \ell$ are given in the two upper rows of the table.
In particular, we have the two possibilities $(a_1,a_2)=(4,1)$ or $(3,2)$.

In the first case, $(a_1,a_2)=(4,1)$, we apply Theorem \ref{thm:bncurves}(i), where the numerical conditions are $g < d^2/8$ and $(d,g) \neq (5,3)$. Conditions
\eqref{eq:min} and
 \eqref{eq:A3} read
\[ d \leq 4 \; \; \mbox{or} \; \; d \geq g+3 \]
and
\[ g \leq 12. \]
Noting that we always have $g < d^2/8$ and $(d,g) \neq (5,3)$ when
$d \geq g+3$, all conditions put together yield
\begin{equation}
  \label{eq:ccc1}
  8g < d^2 \leq 16 \; \; \mbox{or} \; \; g \leq \min\{12,d-3\}.
\end{equation}

In the second case, $(a_1,a_2)=(3,2)$, we apply Theorem \ref{thm:bncurves}(ii), where the numerical conditions are $g=(d^2+3)/12$ or $g < d^2/12$ and $(d,g) \neq (7,4)$. Conditions
\eqref{eq:min} and
 \eqref{eq:A3} read
\[ d \leq 12 \; \; \mbox{or} \; \; g \leq 2d-13 \]
and
\[ g \leq 22. \]
Noting that we always have $g < d^2/12$ and $(d,g) \neq (7,4)$ when
$g \leq 2d-13$, and that the only integer values for $d$ and $g$ satisfying $d \leq 12$ and $g=(d^2+3)/12$ are $(d,g)=(3,1)$ and $(9,7)$, 
all conditions put together yield
\begin{eqnarray}
  \label{eq:ccc2}
   12g < d^2 \leq 144 \;\mbox{with} \; (d,g) \neq (7,4); &  \\ \nonumber \mbox{or} \; \;  (d,g) \in \{(3,1),(9,7)\};   & \mbox{or} \; \;
g \leq \min\{22,2d-13\}.
\end{eqnarray}
Finding, for every $g$, the lowest possible bound for $d$ given by \eqref{eq:ccc1}
and \eqref{eq:ccc2}, yields the numerical conditions in Theorem \ref{result}(a).

\bibliographystyle{amsplain}

\begin{thebibliography}{10}


\bibitem{cl1} H.~Clemens, \textit{Homological equivalence, modulo algebraic
 equivalence, is not finitely generated}, Publ. Math. IHES
 \textbf{58} (1983), 19-38.


\bibitem{ck} H.~Clemens, H.~P.~Kley, 
\textit{Counting curves that move with threefolds},
J. Algebraic Geom. \textbf{9} (2000), 175--200.

\bibitem{dk} I.~Dolgachev, M.~Kapranov,
\textit{Arrangements of hyperplanes and vector bundles on $\PP^n$},
Duke Math. J. {\bf 71} (1993), 633-664. 







\bibitem{EJS} T.~Ekedahl, T.~Johnsen, D.~E.~Sommervoll,
  \textit{Isolated rational curves on K3-fibered Calabi-Yau threefolds},
    Manuscr. Math. \textbf{99} (1999), 111-133.

\bibitem{ev} H.~Esnault, E.~Viehweg, 
  \textit{Lectures on vanishing theorems}, DMV seminar, Band 20, Birkh{\"a}user Verlag, Basel-Boston-Berlin (1992).
    

\bibitem{katz} S.~Katz, \textit{On the finiteness of rational curves on
  quintic  threefolds}, Compos. Math. \textbf{60} (1986), 151-162.




\bibitem{kl1} H.~P.~Kley, \textit{Rigid curves in complete intersection
  Calabi-Yau threefolds}, Compos. Math. \textbf{123} (2000), 185-208.

\bibitem{kl2} H.~P.~Kley, \textit{On the existence of curves in K-trivial
  threefolds}, Preprint, math.AG/9811099 (1998).

\bibitem{kn1} A.~L.~Knutsen, \textit{Smooth curves on projective $K3$
surfaces}, Math. Scand. {\bf 90} (2002), 215-231.

\bibitem{kn2} A.~L.~Knutsen, \textit{Smooth, isolated curves in families of 
Calabi-Yau threefolds in homogeneous spaces}, revised version in preparation.



\bibitem{mori} S.~Mori, \textit{On degrees and genera of curves on smooth
  quartic surfaces in $\PP^3$}, Nagoya Math. J. \textbf{96} (1984), 127-132.


\bibitem{og} K.~Oguiso, \textit{Two remarks on Calabi-Yau Moishezon
  threefolds},
  J. f\"{u}r die Reine und Angew. Math. \textbf{452} (1994), 153-161.

\bibitem{S-D} B.~Saint-Donat, \textit{Projective models of $K-3$
  surfaces}, Amer. J. Math. \textbf{96} (1974), 602-639.



\end{thebibliography}

\end{document}